\documentclass{article}
\usepackage{mathrsfs,amsmath,amssymb,mathtools,amsthm}
\usepackage{tikz-cd}
\usepackage{hyperref}
\usepackage{comment}
\usepackage[capitalise,nameinlink]{cleveref}
\usepackage{dirtytalk}

\DeclareMathOperator{\Mod}{Mod}
\DeclareMathOperator{\Lex}{Lex}
\DeclareMathOperator{\Fix}{Fix}

\DeclareMathOperator{\Fam}{Fam}
\DeclareMathOperator{\PSh}{PSh}

\newcommand{\id}{\textup{id}}
\newcommand{\op}{\textup{op}}

\newcommand{\Set}{\textup{\textbf{Set}}}

\newcommand{\SFr}{\textup{\textbf{S5Fr}}}
\newcommand{\SA}{\textup{\textbf{S5Alg}}}

\newcommand{\Setfs}{\textup{\textbf{Set}}_\textup{s}^\textup{f}}
\newcommand{\fin}{\textup{fin}}
\newcommand{\lf}{\textup{lf}}

\newcommand{\graffe}[1]{\{ #1 \}}

\newcommand{\lra}{\longrightarrow}

\newcommand{\un}{\underline{n}}
\newcommand{\um}{\underline{m}}

\newcommand{\mC}{\mathsf{C}}
\newcommand{\IndC}{\ensuremath{\rm{Ind}\-\mathsf{C}}}

\makeatletter
\newcommand*{\pair}{%
  \mathrel{%
    \mathop{%
      \vcenter{%
        \offinterlineskip
        \halign{%
          \hfil##\hfil\cr
          $\longrightarrow$\cr
          \noalign{\vskip0.1ex}
          $\longrightarrow$\cr
        }%
      }%
    }%
  }%
}
\makeatother

\crefname{theorem}{Theorem}{Theorems}
\crefname{remark}{Remark}{Remarks}
\crefname{corollary}{Corollary}{Corollaries}
\crefname{equation}{}{}

\newtheorem{theorem}{Theorem}[section]
\newtheorem{proposition}{Proposition}[section]
\newtheorem{corollary}{Corollary}[section]

\theoremstyle{definition}
\newtheorem{definition}{Definition}[theorem]
\newtheorem{remark}{Remark}
\newtheorem{example}{Example}

\title{An essentially algebraic glance to Kripke semantics: the 
\textit{S5} case
}
\author{Matteo De Berardinis, Silvio Ghilardi}
\date{October 2025}

\begin{document}

\maketitle

\begin{abstract}
We show that the category of finite \textit{S5}-algebras (dual to finite reflexive, symmetric and transitive Kripke frames) classifies the essentially algebraic theory  whose models are
%can be described as
Kan extensions of faithful actions of the finite symmetric groups.

\end{abstract}

\section{Introduction}
A logic is usually investigated from two basic points of view, namely the syntactic and the semantic points of view. Syntax is usually associated with algebraic structures via Lindenbaum algebras, to the point that algebraic semantics is sometimes referred to as ``syntax in disguise'' %(van Benthem~\cite{johan}) 
or even ``syntax in invariant form'' (Lawvere~\cite{Law}).  On the other hand, since it deals with points, relations, neighborhoods, etc., semantics is usually considered  part of the geometric world. In fact, after the introduction of suitable topological structures in the semantic framework, a perfect duality between (algebraic) syntax and (geometric) semantics  is established. This leads to the well-known dualities between distributive lattices and Priestley spaces, Heyting algebras and Esakia spaces, Boolean algebras and Stone spaces, modal algebras and descriptive frames, etc.

However, there are peculiar situations where the above standard picture 
needs to be somewhat reconsidered.
This is the case of \emph{profinite} algebras. On the one hand, such algebras naturally appear as internal algebraic structures in the category of  Stone spaces~\cite{stone}, so they carry on genuinely geometric features. On the other hand,   monadicity theorems~\cite{our_APAL} still reveal their syntactic character and in fact they are Lindembaum algebras of theories in infinitary languages~\cite{our_mai}. Their duals in the case of modal logic are \emph{locally finite Kripke frames} and p-morphisms; these  are a locally finitely presentable category~\cite{our_APAL}, i.e. they are the category of models of an essentially algebraic theory. According to general facts~\cite{AR}, the syntax of such a theory is represented by the finitely presentable objects, which are nothing but finite Kripke frames  in such  context. Thus, finite Kripke frames are ``syntax in disguise'' or ``invariant syntax'' of an essentially algebraic theory in  this case.

There is a kind of ``syntactic duality'' behind this phenomenon: in fact, when the finitely presented algebras of a variety are closed under equalizers and products, their category (not only the opposite of their  category) classifies an essentially algebraic theory: the two involved theories can be considered  in a sense ``dual'' to each other. This phenomenon arises when the variety is locally finite, but not only in that case (from the results in~\cite{GZ}, it follows that Heyting algebras are a non locally finite example). Identifying the `dual theory' of the equational theory axiomatizing a variety (when such theory exists) might be tricky, beyond  easy cases like those listed in \cref{sec:ea} below. The case study presented in this paper (monadic Boolean algebras, also known as \textit{S5} modal algebras) is a first non trivial example in this sense.

The paper is structured as follows: in \cref{sec:ea} we recall the ingredients of the Gabriel-Ulmer theory that are needed to understand the subsequent sections.  \cref{sec:ma} introduces a digression on Ind-completions, modal algebras and locally finite Kripke frames. \cref{sec:families} reports some basic facts about the finitary version of the  construction of the category of families.
In \cref{sec:actions} we introduce symmetric groups actions and their canonical liftings. \cref{sec:classified} contains our axiomatization result and \cref{sec:conclusions} reports some concluding remarks.

\section{Essentially algebraic theories}\label{sec:ea}

Essentially algebraic theories arise in the context of Gabriel-Ulmer duality theory~\cite{GU} (see~\cite{MP} for a synthetic exposition). We give here essential information concerning some syntactic and categorical aspects, referring to~\cite[Chapter D1]{Elephant} for more 
(alternative approaches, using languages with partial functions,
can be found in~\cite{AR,PV}). 
The signature of an essentially algebraic theory is a customary first-order multi-sorted signature with equality; without loss of generality, one may take the signature to be purely relational or purely functional, as shown in~\cite[Lemma D1.4.9]{Elephant}
(in this paper, we limit ourselves to purely functional signatures, unless otherwise stated).

In an essentially algebraic theory we can build formulae using equalities, conjunctions and existential quantifiers, but an existential quantifier like $\exists y B(\underline{x},y)$  can be introduced   only under the proviso that the theory proves the uniqueness condition 
$$
B(\underline{x},y_1) \wedge B(\underline{x},y_2) \to y_1 = y_2~~.
$$ 
As axioms, an essentially algebraic theory comprises sequents of the form $A\to B$,
where $A$, $B$ are formulae built up according to the above instructions.
Thus, the introduction of the set of formulae and of the set of axioms needs a mutual recursion~\cite[\S D1.3.4, p.833]{Elephant}. 
%this is the approach we follow, taken e.g. from~\cite[D1.3.4, p.833]{Elephant}.\footnote{Alternative definitions are possible, see Adamek-Rosicky 1994 or Palmgren-Vickers 2007.}

The categories of models of essentially algebraic theories are equivalent to the categories of \Set-valued lex functors having a small lex category as domain (the word `lex' refers to finite limits and their preservation).
In fact, given a small lex category  $\mathsf{C}$, we can associate with it an essentially algebraic theory $T_\mathsf{C}$ as follows: we have one sort for every object $A$ of 
$\mathsf{C}$ and one unary function symbol $f$ of domain sort $A$ and codomain sort $B$ for every arrow $f \colon A \lra B$ in $\mathsf{C}$. We can express functoriality with axioms of the kind
$f(g(x))= (f \circ g)(x)$; to express preservation of finite limits, we can use essentially algebraic axioms. For instance, if $e \colon E \lra A$ is the equalizer of $f,g \colon A \pair B$, we first state injectivity of $e$ via
$$
e(x_1)=e(x_2) \to x_1=x_2
$$
and then we write the equalizer condition as
$$
f(y)=g(y) \to \exists x(e(x)=y)~~.
$$
%matteo
Given a diagram of pullback
\[\begin{tikzcd}
	P & {A_2} \\
	{A_1} & B
	\arrow["{p_2}", from=1-1, to=1-2]
	\arrow["{p_1}"', from=1-1, to=2-1]
	\arrow["{f_2}", from=1-2, to=2-2]
	\arrow["{f_1}"', from=2-1, to=2-2]
\end{tikzcd}\]
its preservation is expressed by the joint injectivity axiom
$$
p_1(x) = p_1(x') \wedge p_2(x) = p_2(x') \to x = x'
$$
together with the pullback condition
$$
f_1(y_1) = f_2(y_2) \to \exists x\,(p_1(x) = y_1 \wedge p_2(x) = y_2)~~.
$$
For products there is a simpler solution: for every pair of objects $A,B$, one introduces a binary pairing function symbol $\langle -,-\rangle$
of domain types $A\,B$ and codomain type $A\times B$ with the axioms
$$
\textit{pr}_1(\langle x,y\rangle)=x, \qquad  \textit{pr}_2(\langle x,y\rangle)=y, \qquad
\langle \textit{pr}_1(z),\textit{pr}_2(z)\rangle=z,
$$
where $\textit{pr}_1, \textit{pr}_2$ are the function symbols corresponding to the product projections $\textit{pr}_1 \colon A\times B \lra A$ and $\textit{pr}_2 \colon A\times B \lra B$
in $\mathsf{C}$.
%matteo
We can also express preservation of the terminal object $1$ of $\mathsf{C}$, by adding a constant $*$ of type $1$ and writing the axiom $x = *$, with $x$ variable of type $1$.

It is then clear that the models of $T_\mathsf{C}$ and the homomorphisms among them\footnote{Here by a homomorphism we mean just a sort-indexed family of functions preserving the operations.} are a category equivalent to the category of lex functors $T_\mathsf{C}\lra \Set$  and natural transformations
(indeed, this statement can be generalized by taking internal models in any other lex category in place of $\Set$).

Vice versa, one can produce out of an essentially algebraic theory $T$, a small lex category 
$\mathsf{C}_T$ so that the category of lex functors from $\mathsf{C}_T$ to $\Set$  is equivalent to the category of models of $T$:
the construction (which is not relevant for the purposes of this paper) is explained in detail in~\cite[D1.4]{Elephant}.
The two passages, from $T$ to 
$\mathsf{C}_T$ and from $\mathsf{C}$ to $T_\mathsf{C}$ are inverse to each other, modulo appropriate notions of equivalence.

For an essentially algebraic theory $T$, we let $\Mod(T)$ to be the category of set-theoretic models of $T$ and of homomorphisms among them: such a category is called \emph{locally finitely presentable} (an intrinsic  definition of locally finitely presentable category is also possible, see~\cite{AR}).
$\mathsf{C}_T$ (better, a lex category equivalent to it)  can be directly recovered from
$\Mod(T)$: indeed,  $\mathsf{C}_T$ is equivalent to the opposite of the full subcategory of $\Mod(T)$ formed by the \emph{finitely presented} models. Finitely presented models can be defined in terms of generators and relations as in the standard algebraic context~\cite{PV}, but there is also an equivalent categorical characterization: $M$ is finitely presented iff the representable functor $\Mod(T)[M,-]$ preserves filtered colimits.
Thus, for instance, if $T$ is the theory of groups, then we have that
\begin{equation}\label{eq:main}
 \Mod(T) ~\simeq~ \Lex(\mathsf{C}_T^\op, \Set)   
\end{equation}
where $\mathsf{C}_T$ is the category of finitely presented groups.
The equivalence of categories \cref{eq:main} holds for every essentially algebraic theory $T$ and constitutes the starting point for the investigations of this paper.

In some cases, the finitely presented models of an essentially algebraic theory $T$ are closed under products and equalizers (this is always the case when $\Mod(T)$ is a locally finite variety, i.e. an equational class of algebras, whose finitely generated members are finite). In this case, \emph{both $\mathsf{C}_T$ and $\mathsf{C}_T^\op$ classify an essentially algebraic theory} in the sense that there is a theory $T^d$ (called \emph{dual} of $T$) such that we have  $\Mod(T^d) ~\simeq~ \Lex(\mathsf{C}_T, \Set)$. Here you are some examples.

\begin{example} According to \cref{eq:1},
the category of finite sets classifies the theory of Boolean algebras, because the opposite of the category of finitely presented (= finite) Boolean algebras is equivalent to the category of finite sets. If we take the dual of the category of finite sets (namely, the category of finite Boolean algebras) we realize that, for the same reason, this is also a classifying category of an essentially algebraic theory: in fact, it classifies the pure equality theory (in a one-sorted language). Thus, the dual of the theory of Boolean algebras is the pure equality theory.
\end{example}

\begin{example}
The theory of join-semilattices (with zero) is dual to the theory of meet-semilattices (with unit): thus, such a theory is in fact self-dual. This is because finitely presented join semilattices are finite and hence complete, so they are lattices. Taking posetal right adjoints, it is clear that we can convert (up to natural bijections)  
 join semilattices morphisms to  meet semilattices morphisms and vice versa.
\end{example}

\begin{example}
The variety of Heyting algebras generated by the 3-element chain has as dual the theory of the actions of the multiplicative monoid of $\mathbb{Z}_2$ (few more examples like this can be found in~\cite{Copr}).
\end{example}

\begin{example}
 This example shows that the dual of an algebraic (=equational) theory need not be algebraic (it is essentially algebraic). Consider  distributive lattices (with zero and unit); finitely presented distributive lattices are finite and hence, by Priestley duality~\cite{P}, they are duals to finite posets. This shows that  that the dual of the theory of distributive lattices is the theory of posets, because finite posets are the finitely presented objects in the category of posets. The theory of posets  can be seen as an essentially algebraic theory in a functional signature as follows, following the proof of
 ~\cite[Lemma D1.4.9]{Elephant}. One takes a two-sorted signature with sorts $X,P$ and two unary function symbols $p_1,p_2 \colon X\pair P$. 
 We have axioms
 $$
 p_1(x)= p_1(x') \wedge p_2(x)= p_2(x') \to x=x'~~
 $$
 allowing us to define $y_1 \leq y_2$ as $\exists x\, (p_1(x)=y_1 \wedge p_2(x)=y_2)$ and to specify reflexivity, transitivity and antisymmetry conditions in the essentially algebraic language.
\end{example}

\subsection{Our case study}

Given a locally finite variety $\mathbf{V}$ (in a finitary language), whereas we know from \cref{eq:main} that the opposite of the category 
%$\mathbf{V}_{\fin}^{op}$ 
of finite $\mathbf{V}$-algebras classifies the equational theory whose models are the algebras in $\mathbf{V}$,  it is not clear what finite $\mathbf{V}$-algebras themselves classify: the dual theory is in general difficult to identify, beyond  simple examples like the above ones. In this paper, we propose a relatively simple (but non trivial) case study: \textit{S5}-algebras.

A \emph{\textit{S5}-algebra} is a Boolean algebra endowed with a unary operator $\Box$ such that:
$$\Box \top = \top, \quad \Box(x \wedge y) = \Box x \wedge \Box y, \quad \Box x \leq x, \quad x \leq \Box \Diamond x$$
where $\Diamond x$ abbreviates $\neg \Box \neg x$. 
\textit{S5}-algebras are also called \emph{monadic Boolean algebras}~\cite{halmos} because these algebras are `Lindenbaum algebras' of  one-variable classical first order logic (the $\Box$ operator is read as the universal quantifier in such interpretation).

We let $\SA$ be the category of \textit{S5}-algebras and homomorphisms (this  is a variety)  and
%\textit{S5}-algebras are a variety, 
we call $T_{\textit{S5}}$ the corresponding equational theory (in the language of Boolean algebras enriched with the unary function symbol $\Box$).

Typical \textit{S5}-algebras come from \textit{S5}-Kripke frames.
An \emph{S5-Kripke frame} is a pair $(W,R)$, where $W$ is a set and $R$ is  an equivalence relation; a morphism between 
\textit{S5}-Kripke frames (a \emph{p-morphism} in the modal logic terminology)
$$
f \colon (W,R) \lra (W', R')
$$
is a function mapping $R$-equivalence classes \emph{onto} $R'$-equivalence classes. 
We let $\SFr$ be the category of \textit{S5}-Kripke frames and morphisms.
An \textit{S5}-Kripke frame is turned into an \textit{S5}-algebra by endowing the powerset Boolean algebra 
$\mathcal{P}(W)$ with the operator
$$
\Box_R \colon \mathcal{P}(W)\lra \mathcal{P}(W)
$$
defined by 
$$
\Box_R(S)~=~\{w\in W\mid \forall v\, (wRv\to v\in S)\}
$$
for every $S\subseteq W$.

The following well-known result can be found in  modal logic textbooks~\cite{misha}.
\begin{theorem}
\textit{S5}-algebras are a locally finite variety. 
The category $\SA_\fin$ of finite \textit{S5}-algebras is dual to the category $\SFr_\fin$ of finite \textit{S5}-Kripke frames.    
\end{theorem}

One might be tempted to say that the dual theory of
$T_\textit{S5}$ 
%the theory of \textit{S5}-algebras 
(namely the essentially algebraic theory classified by finite \textit{S5}-algebras as a lex category)
is the theory of a set endowed with an equivalence relation. This is not correct, because such a theory is classified by the opposite of the category having as objects finite sets endowed with an equivalence relation, but as arrows the functions preserving the equivalence relation, i.e.\ mapping equivalent pairs to equivalent pairs (clearly not all such functions map equivalence classes onto equivalence classes). In other words, the key point is that not all relation-preserving maps are p-morphisms. 
In fact, the essentially algebraic theory classified by finite \textit{S5}-algebras requires a much more complex specification.

\section{Ind-Completions and Local Finiteness}\label{sec:ma}

In this Section, we report some material from~\cite{our_APAL}. Using it, we shall be able to give a description  of  \emph{the category of models} of the dual theory of $T_\textit{S5}$. This is interesting information by itself, but unfortunately such a 
description will not lead to an explicit \emph{logical presentation} of $T_\textit{S5}^d$: the latter requires the design of a  set of axioms in a suitable first-order language  and 
will be supplied only in \cref{sec:classified} (the present section can be skipped in a first reading).

Let $\mathbf{V}$ be a locally finite variety.
 % (in a finitary language)
 From general facts, one can prove that $\Lex(\mathbf{V}_{\fin},\Set)$ is the opposite of the category of profinite $\mathbf{V}$-algebras (the latter being the Pro-Completion of $\mathbf{V}_{\fin}$, the category of finite $\mathbf{V}$-algebras). 
 The notion of Pro-completion is dual to the notion of Ind-completion: we are going to recall the latter, taking some information from~\cite[Chapter VI]{stone}.

 The Ind-Completion of a category $\mC$ formally adds filtered colimits to $\mC$. One way of introducing it is as a  full subcategory of the presheaf category:

\begin{definition} For a small category $\mC$, we let  $\IndC$ be
 the full subcategory of $\Set^{\mC^\op}$ given by those functors that are filtered colimits of representable functors.
\end{definition}

We are only interested in the case where $\mC$ has finite colimits;
the following theorem
collects various statements from~\cite[Chapter VI, \S 1]{stone} and
summarizes the relevant properties of \IndC\ in
that hypothesis.
\begin{theorem}\label{thm:indchar}
Let $\mC$ be a small category with finite colimits. Then
\begin{itemize}
    \item[{\rm (i)}] \IndC\ is equivalent to $\Lex(\mC^\op,\Set)$, i.e.\ to the full subcategory of $\Set^{\mC^\op}$ given by the contravariant functors from $\mC$ to \Set\ turning finite colimits into finite limits;
    \item[{\rm (ii)}] \IndC\ has all filtered colimits and the embedding $\IndC \longrightarrow \Set^{{\mC}^\op}$ preserves them;
    \item[{\rm (iii)}] \IndC\ has  finite colimits and  the Yoneda embedding (restricted in the co\-do\-main)  $Y \colon \mC \longrightarrow \IndC$ preserves them;
    \item[{\rm (iv)}] \IndC\  is complete, and the embedding $\IndC \longrightarrow \Set^{{\mC}^\op}$ preserves all small limits.
   % \item[{\rm (iv)}] the embedding $Y : \mC \longrightarrow \IndC$ preserves all limits which exist in $\mC$.
\end{itemize}
\end{theorem}

Notice that, as a consequence of (ii) and (iii), we have that \IndC\ is co-complete;
actually it is a co-completion of $\mC$, in the following sense.
\begin{definition}
Given a small category $\mC$, a \textit{cocompletion} of $\mC$ is a full embedding $F \colon \mC \longrightarrow \mathsf{D}$ into a cocomplete category $\mathsf{D}$ such that every object of $\mathsf{D}$ is a colimit of objects in the image of $F$.
\end{definition}

% We recall the notion of finite presentability.

%\begin{definition}
%An object $X$ of a (locally small) category with filtered colimits $\mathcal{D}$ is said to be \textit{finitely-presentable} (in $\mathcal{D}$) if the functor $\text{Hom}_{\mathcal{D}}(X,-) : \mathcal{D} \rightarrow \Set$ preserves filtered colimits.
%\end{definition}

The following theorem
(taken literally from~\cite[Chapter VI, Theorem 1.8]{stone}) characterizes \IndC\ up to equivalence as that cocompletion of $\mC$ for which the embedding functor preserves finite colimits and sends the objects of $\mC$ to finitely-presentable objects.

\begin{theorem}\label{teo:charindcompl}
Let $\mC$ be a small category with finite colimits, and $Z \colon \mC \longrightarrow \mathsf{D}$ a full embedding of $\mC$ in a cocomplete category $\mathsf{D}$. Then
\begin{itemize}
    \item[{\rm (i)}] if the objects in the image of $Z$ are finitely-presentable in $\mathsf{D}$, $Z$ extends to a full embedding $\hat Z \colon \IndC \longrightarrow \mathsf{D}$;
    \item[{\rm (ii)}] if in addition $Z \colon \mC \longrightarrow \mathsf{D}$ is a cocompletion of $\mC$ and $Z$ preserves finite colimits, $\hat Z$ is an equivalence.
\end{itemize}
\end{theorem}

Thus,
given a variety  $\mathbf{V}$, one can apply \cref{teo:charindcompl} to show that \emph{$\mathbf{V}$ is the Ind-Completion of the category of finitely
presented $\mathbf{V}$-algebras}~\cite[Chapter VI, Corollary 2.2]{stone}. 
%
%We now take some results from~\cite{our_APAL} and show how the notion of Ind-Completion can be applied in dual semantic contexts.
%In the area of modal logic, there are nice dualities between finite frames and finite algebras: the idea is to exploit these finite dualities in order to better understand the essentially algebraic theories classified by finitely presented algebras in locally finite varieties. The outcome of such a research program could be a better understanding of locally finite Kripke frames (namely the Ind-completion of finite Kripke frames) as a category equivalent to a class of essentially algebraic structures. In these notes, we analyze a rather simple but non trivial case, where $\mathbf{V}$ is the variety of modal 
\cref{teo:charindcompl} can also be used in order to obtain the next result.
\begin{theorem}\label{thm:fromAPAL}
An \textit{S5}-Kripke frame $(W, R)$ is said to be \emph{locally finite} if and only if the $R$-equivalence classes all have finite cardinality; let $\SFr_\lf$  be the full subcategory of $\SFr$ formed by locally finite \textit{S5}-Kripke frames. Then $\SFr_\lf$ is the Ind-Completion of $\SFr_\fin$ and it is equivalent to $\Lex(\SFr_\fin^\op,\Set)$.
\end{theorem}

\begin{proof}
The proof of this statement (in the more general case of the Kripke frames for the modal system $K$) is given in~\cite[Sections 3-4]{our_APAL}. We report it here for the sake of completeness.

First, observe that $\SFr_\fin$ is equivalent to its full subcategory of the sets $\underline{n} = \{1,\dots,n\}$ endowed with some equivalence relation, which is a small category.

We prove that $\SFr$ is cocomplete and that the forgetful functor $\SFr \longrightarrow \Set$ preserves all colimits. Moreover, we prove that the embeddings $\SFr_\lf \longrightarrow \SFr$ and $\SFr_\fin \longrightarrow \SFr$ create  all colimits and all finite colimits, respectively. The initial object in $\SFr$ is the empty set endowed with the empty relation. The coproduct in $\SFr$ of a family $\{(W_i,R_i)\}_{i \in I}$ of equivalence relations is the disjoint union $(\coprod_{i \in I} W_i,\coprod_{i \in I} R_i)$. The inclusions $\iota_j \colon P_j \longrightarrow \coprod_{i \in I} W_i$ define p-morphisms, for all $j \in I$ and the universal property of the coproduct is easily checked. Moreover, if all the $W_i$ are (locally) finite, then so it is their coproduct. The coequalizer in $\SFr$ of a pair of parallel p-morphisms $f,g \colon (W,R) \pair (V,S)$ is given by the quotient of $V$ via the equivalence relation $\approx$ generated by the set of pairs $\{(f(w),g(w))\ \vert\ w \in W\}$, endowed with the equivalence relation ${S}/{\approx}$ defined as follows: for $[x], [y] \in {V}/{\approx}$,
\begin{align*}
[x] ({S}/{\approx}) [y] \iff \exists x', y' \in V \text{ s.t.\ } x \approx x',\ x' S y' \text{ and } y' \approx y
\end{align*}
It is possible to show (see \cite[Proposition 3.9]{our_APAL}) that, with such definitions, the projection $q \colon V \longrightarrow {V}/{\approx}$ defines a p-morphism; obviously, $q \circ f = q \circ g$. Moreover, given any other p-morphism $p \colon (V,S) \longrightarrow (U,T)$ such that $p \circ f = p \circ g$, then there exists a unique function $h \colon {V}/{\approx} \longrightarrow U$ such that $h \circ q = p$. The fact that both $q$ and $p$ are p-morphisms forces $h$ to be so, too. To conclude, observe that $({V}/{\approx},{S}/{\approx})$ is (locally) finite if $(V,S)$ is so.

We now prove that the objects of $\SFr_\fin$ are finitely-presentable in $\SFr_\lf$. Consider a colimiting cocone $(\iota_i \colon (W_i,R_i) \longrightarrow (W,R)\ \vert\ i \in I)$ in $\SFr_\lf$ over some filtered diagram $I \longrightarrow \SFr_\fin$ of finite S5-frames $(W_i,R_i)$. We have to prove that, for each $(V,S) \in \SFr_\fin$, any p-morphism $f \colon (V,S) \longrightarrow (W,R)$ factorizes in $\SFr_\lf$ as $(V,S) \buildrel{f_i}\over\longrightarrow (W_i,R_i) \buildrel{\iota_i}\over\longrightarrow (W,R)$ in a unique way up to equivalence ($f_i \colon (V,S) \longrightarrow (W_i,R_i)$ and $f_j \colon (V,S) \longrightarrow (W_j,R_j)$ being equivalent if there exist $i \longrightarrow k$ and $j \longrightarrow k$ in $I$ such that the induced diagram
\[\begin{tikzcd}[ampersand replacement=\&]
	{(V,S)} \& {(W_i,R_i)} \\
	{(W_j,R_j)} \& {(W_k,R_k)}
	\arrow["{f_i}", from=1-1, to=1-2]
	\arrow["{f_j}"', from=1-1, to=2-1]
	\arrow[from=1-2, to=2-2]
	\arrow[from=2-1, to=2-2]
\end{tikzcd}\]
commutes, see \cite[Chapter VI]{stone}). Since the forgetful functor $\SFr_\lf \longrightarrow \Set$ preserves colimits and finite sets are finitely-presentable in $\Set$, the existence and essentially uniqueness of a function $f_i \colon V \longrightarrow W_i$ realizing the desired factorization is guaranteed. In general, $f_i$ is just a function. We prove that $f_i$ is a p-morphism for a suitable choice of the index $i$. We first observe that, since $\SFr_\lf \longrightarrow \Set$ preserves colimits, $W$ is isomorphic, as a set, to a quotient of $\coprod_{i \in I} W_i$; moreover, for $i, j \in I$, $x \in W_i$ and $x' \in W_j$, we have that $\iota_i(x) = \iota_j(x')$ if and only if there exist $i \longrightarrow k$ and $j \longrightarrow k$ in $I$ such that the induced $W_i \longrightarrow W_k$ and $W_j \longrightarrow W_k$ send $x$ and $x'$ to the same element. If we have $v, v' \in V$ such that $v S v'$, then ($f$ is a p-morphism) $\iota_i(f_i(v)) = f(v)Rf(v') = \iota_i(f_i(v'))$. Since $\iota_i$ is a p-morphism, there exists $x' \in W_i$ such that $f_i(v)R_ix'$ and $\iota_i(x') = \iota_i(f_i(v'))$. Using the previous remark and the fact that any two parallel arrows $i \pair j$ in $I$ can be coequalized by an arrow $j \longrightarrow k$ (by the properties of a filtered category, see \cite[Chapter VI]{stone}), we obtain a morphism $i \longrightarrow k$ in $I$ such that the induced $W_i \longrightarrow W_k$ sends $x'$ and $f_i(v')$ to 
the same element; this means that, up to composition with $W_i \longrightarrow W_k$, we can assume that $f_i(v) R_i f_i(v')$. Moreover, if a certain pair $v S v'$ is preserved by $f_i$, then the composition of $f_i$ with $W_i \longrightarrow W_k$ (which is a p-morphism) preserves it, too. The set $W$ being finite, we can assume that $f_i$ maps equivalence relations to equivalence relations. Let us now prove that equivalence relations are mapped onto equivalence relations. If we have $f_i(v) R_i x'$ in $W_i$, then ($\iota_i$ is a p-morphism) $f(v) = \iota_i(f_i(v)) R \iota_i(x')$, hence ($f$ is a p-morphism) there exists $v' \in V$ such that $v S v'$ and $\iota_i(f_i(v')) = f(v') = \iota_i(x')$. As before, if we compose with $W(d) \colon W_i \longrightarrow W_k$ (induced by a suitable $d \colon i \longrightarrow k$ in $I$), we have that $W(d)(f_i(v')) = W(d)(x')$. We proved that
\begin{center}
(*) \say{for all $v \in V$, for all $x' \in W_i$ such that $f_i(v)R_ix'$, there are $v' \in V$ and $d \colon i \longrightarrow k$ in $I$ s.t. $v S v'$ and $W(d)(f_i(v')) = W(d)(x')$}.
\end{center}
Since $V$ is finite and $W_i$ is locally finite, there are finitely many such pairs $(v,x')$ and, consequently, $I$ being filtered, we can take the same $d$ for all such pairs. The composition $W(d) \circ f_i$ now fits our purposes: if we have $W(d)(f_i(v))R_ky'$ for some $y' \in W_k$, then (as $W(d)$ is a p-morphism) there is $x' \in W_i$ such that $f_i(v)R_ix'$ and 
$W(d)(x') = y'$: applying (*) to the pair $(v,x')$ we get $v' \in V$ such that $v S v'$ and $W(d)(f_i(v')) = W(d)(x') = y'$, proving that $W(d) \circ f_i$ maps equivalence classes onto equivalence classes (in other words, $f_i \colon V \longrightarrow W_i$ itself does so, for a suitable choice of the index $i$).

To conclude, it is sufficient to observe that every object of $\SFr_\lf$ is a colimit of objects in $\SFr_\fin$, since it can be written as the disjoint union, i.e.\ the coproduct, of its equivalence classes, which are finite by definition.
\end{proof}

Thus the category of models of the dual theory of $T_\textit{S5}$ can be represented in two equivalent ways: as $\SFr_\lf$ or as 
$\Lex(\SFr_\fin^\op,\Set) \simeq \Lex(\SA_\fin,\Set)$. Both ways can be useful: for instance, limits are easy to describe using 
the latter representation
%$Lex(\SFr_\fin^{op},\Set)\simeq Lex(\SA_\fin,\Set)$ 
(but colimits are hard, as they go through the reflection functor from presheaves). On the other hand, colimits are easy to describe in 
$\SFr_\lf$ (because they are preserved by the forgetful functor to $\Set$ as shown in the above proof),
%~\cite{our_APAL}),
but products are tricky (see~\cite{Copr}).

\section{The finite families category}\label{sec:families}
%We know that $\SFr_\lf$ is equivalent to the category $\Lex(\SFr_\fin^\op,\Set)$ of lex functors  $\SA_\fin \longrightarrow \Set$. 
%
In this section, we give another presentation of 
$\Lex(\SFr_\fin^\op,\Set)$
%$\SFr_\lf$ 
as a category of presheaves, this time considering (the dual of) a smaller class of algebras as the domain site.

%silvio 28/10 riscrivo lievemente questo paragrafo perchè avevamo detto che il paragrafo precedente poteva essere saltato in prima lettura
We denote by $\Setfs$ the category having as objects the sets $\underline{n} = \{1,\dots,n\}$ and as maps the surjective functions among them (this is equivalent to the category of non-empty finite sets and surjections). There exist a full embedding $\Setfs \longrightarrow \SFr_\fin$, 
%\longrightarrow \SFr_\lf$, 
sending $\underline{n}$ to the $n$-elements cluster: this is the finite Kripke frame $(\underline{n}, U)$, where $U=\underline{n}\times \underline{n}$ is the total relation.\footnote{This embedding represent  $\Setfs$ as the dual of the category of finite subdirectly irreducible \textit{S5}-algebras.} From another point of view, the full embedding $\Setfs \longrightarrow \SFr_\fin$ presents $\SFr_\fin$ as the category obtained by freely adjoining  finite coproducts to $\Setfs$. In a similar way, the full embedding $\Setfs \longrightarrow \SFr_\lf$ into the category of locally finite \textit{S5}-frames we met in the previous section 
% $\Setfs \longrightarrow \SFr_\lf$ (resp. ) 
presents $\SFr_\lf$ as the category obtained by freely adjoining all coproducts  to $\Setfs$. 
 
 These are instances of the general construction of the \emph{category of (finite) families} that we are going to recall.  Given a small category $\mathsf{C}$, we can embed it into the category $\Fam(\mathsf{C})$ ($\Fam_\fin(\mathsf{C})$) of (finite) families of $\mathsf{C}$. Namely, an object of $\Fam(\mathsf{C})$ ($\Fam_\fin(\mathsf{C})$) is a gadget of the form $X_I = (I, (X_i)_{i \in I})$, with $I$ a (finite) set and $X_i$ an object of $\mathsf{C}$ for each $i \in I$. An arrow $X_I \longrightarrow Y_J$ is given by a pair $(f,\phi)$, where $f \colon I \longrightarrow J$ is a function and $\phi$ is a family of morphisms $\phi_i \colon X_i \longrightarrow Y_{f(i)}$ in $\mathsf{C}$.
\begin{proposition}\label{prop:free}
If $\mathsf{E}$ has all (finite) coproducts, then any functor $F \colon \mathsf{C} \longrightarrow \mathsf{E}$ extends (essentially uniquely) to a functor $\hat{F} \colon \Fam(\mathsf{C}) \longrightarrow \mathsf{E}$ ($\hat{F} \colon \Fam_\fin(\mathsf{C}) \longrightarrow \mathsf{E}$) that preserves all (finite) coproducts.
\end{proposition}
\begin{proof} The proof is via left Kan extension~\cite[Prop.2.1 and Rem.2.2]{hu-tholen}, we give here a direct elementary description.
Define $\hat{F} (X_I)$ to be the $\mathsf{E}$-coproduct $\coprod_{i \in I} F(X_i)$. The universal property of the coproduct allows us to extend such assignment to a functor $\Fam(\mathsf{C}) \longrightarrow \mathsf{E}$ ($\Fam_\fin(\mathsf{C}) \longrightarrow \mathsf{E}$), in such a way that the restriction to $\mathsf{C}$ is isomorphic to $F$: given $(f,\phi) \colon X_I \longrightarrow Y_J$ in $\Fam(\mathsf{C})$ ($\Fam_\fin(\mathsf{C})$), if we denote by $\iota_j \colon F(Y_j) \longrightarrow \coprod_{j \in J} F(Y_j) = \hat{F}(Y_J)$ the coprojections, then the family of morphisms $\iota_{f(i)} \circ F(\phi_i) \colon F(X_i) \longrightarrow F(Y_{f(i)}) \longrightarrow \hat{F}(Y_J)$, varying $i \in I$, induces a unique morphism $\hat{F}(f,\phi) \colon \hat{F}(X_I) \longrightarrow \hat{F}(Y_J)$ in $\mathsf{E}$. Obviously, $\hat{F}$ preserves (finite) coproducts, and its definition is the only one for which (finite) coproducts are preserved.
\end{proof}

In general, the embedding $\mathsf{C} \longrightarrow \Fam_\fin(\mathsf{C})$ does not preserve existing finite coproducts; however, other existing colimits may be preserved.
%matteo
In particular, we will consider coequalizers and pushouts (and their duals, equalizers and pullbacks). We call \emph{finite connected colimits} (\emph{finite connected limits}) the colimits (limits) generated by them.\footnote{By (co)limits \textit{generated} by some given primitive (co)limits, we mean those (co)limits that can be obtained after applying the primitive (co)limits a finite number of times. Existence and preservation properties of such (co)limits can then be checked on the primitive ones.}
\begin{proposition}\label{prop:pres}
If $\mathsf{C}$ has finite connected colimits, then so does $\Fam_\fin(\mathsf{C})$, and the canonical embedding $\mathsf{C} \longrightarrow \Fam_\fin(\mathsf{C})$ preserves them. Moreover, if $\mathsf{E}$ is finitely cocomplete and $F \colon \mathsf{C} \longrightarrow \mathsf{E}$ preserves finite connected colimits, then the extension $\hat{F} \colon \Fam_\fin(\mathsf{C}) \longrightarrow \mathsf{E}$ preserves all finite colimits.
\end{proposition}
\begin{proof}
We start by describing coequalizers and pushouts in $\Fam_\fin(\mathsf{C})$. For coequalizers, consider a diagram $(f,\phi), (g,\gamma) \colon X_I \pair Y_J$ in $\Fam_\fin(\mathsf{C})$.

We first show that, without loss of generality, we can assume $f = g$. For, if there exists $i \in I$ such that $f(i) \neq g(i)$, define the set $J'$ to be the quotient of $J$ over the equivalence relation generated by the pair $(f(i),g(i))$; call $p \colon J \longrightarrow J'$ the projection map. Then, for each $j \in J$, the morphism $\pi_j \colon Y_j \longrightarrow Y'_{p(j)}$ is defined to be:
\begin{enumerate}
    \item the identity, if $j \notin \{f(i),g(i)\}$;
    \item one of the morphisms obtained by pushing out in $\mathsf{C}$ the diagram
    \[\begin{tikzcd}
	{Y_{g(i)}} & {X_i} & {Y_{f(i)},}
	\arrow["{\gamma_i}"', from=1-2, to=1-1]
	\arrow["{\phi_i}", from=1-2, to=1-3]
    \end{tikzcd}\]
    otherwise.
\end{enumerate}
We defined a morphism $(p,\pi) \colon Y_J \longrightarrow Y'_{J'}$ in $\Fam_\fin(\mathsf{C})$, with the following property: for any morphism $(d,\delta) \colon Y_J \longrightarrow T_L$ such that $(d,\delta) \circ (f,\phi) = (d,\delta) \circ (g,\gamma)$, there exists a unique morphism $(e,\epsilon) \colon Y'_{J'} \longrightarrow T_L$ such that $(e,\epsilon) \circ (p,\pi) = (d,\delta)$ (it means that computing the coequalizer in $\Fam_\fin(\mathsf{C})$ of $(f,\phi)$ and $(g,\gamma)$ it's the same as computing the coequalizer of their compositions with $(p,\pi)$). This allows us to replace $(f,\phi)$ and $(g,\gamma)$ with $(p,\pi) \circ (f,\phi)$ and $(p,\pi) \circ (g,\gamma)$, so that the number of $j \neq j' \in J$ such that $j = f(i)$ and $j' = g(i)$ for some $i \in I$ strictly decreases. This process terminates, the set $J$ being finite.

Now, the coequalizer of $(f,\phi), (f,\gamma) \colon X_I \pair Y_J$ in $\Fam_\fin(\mathsf{C})$ is given by $(\id_J,\eta) \colon Y_J \longrightarrow Z_J$, where $\eta_j \colon Y_j \longrightarrow Z_j$ is:
\begin{enumerate}
    \item the identity, if $j \notin f[I]$;
    \item the iterated $\mathsf{C}$-coequalizer of all the $\phi_i, \gamma_i \colon X_i \pair Y_{f(i)}$, for each $i \in f^{-1}(j)$, otherwise (this process terminates, the set $I$ being finite).
\end{enumerate}
Obviously, coequalizers in $\mathsf{C}$ are preserved by $\iota$. Similarly, pushouts in $\Fam_\fin(\mathsf{C})$ can be computed by taking iterated pushouts in $\mathsf{C}$, and pushouts in $\mathsf{C}$ are preserved by $\iota$.

Finally, we need to check that $\hat{F}$ of \cref{prop:free} preserves coequalizers (together with the preservation of finite coproducts, we have the preservation of all finite colimits). Again, given $(f,\phi), (g,\gamma) \colon X_I \pair Y_J$ in $\Fam_\fin(\mathsf{C})$, we can assume $f = g$. Let us 
%unravel the process in
consider the construction from the first part of the proof; applying $\hat{F}$, we obtain a diagram in $\mathsf{E}$
\[\begin{tikzcd}
	{\hat{F}(X_I)} & {\hat{F}(Y_J)} & {\hat{F}(Y'_{J'})}
	\arrow["{\hat{F}(f,\phi)}", shift left, from=1-1, to=1-2]
	\arrow["{\hat{F}(g,\gamma)}"', shift right, from=1-1, to=1-2]
	\arrow["{\hat{F}(p,\pi)}", from=1-2, to=1-3]
\end{tikzcd}\]
We prove the following: for any morphism $a \colon \hat{F}(Y_J) \longrightarrow E$ in $\mathsf{E}$ such that $a \circ \hat{F}(f,\phi) = a \circ \hat{F}(g,\gamma)$, there exists a unique morphism $b \colon \hat{F}(Y'_{J'}) \longrightarrow E$ such that $b \circ \hat{F}(p,\pi) = a$ (it means that computing the coequalizer in $\mathsf{E}$ of $\hat{F}(f,\phi)$ and $\hat{F}(g,\gamma)$ it's the same as computing the coequalizer of their compositions with $\hat{F}(p,\pi)$). Take the components $a_j \coloneqq a \circ \iota_j \colon F(Y_j) \longrightarrow \hat{F}(Y_J) \longrightarrow E$:
\begin{enumerate}
    \item If $j \notin \{f(i),g(i)\}$, then $F(\pi_j) \colon F(Y_j) \longrightarrow F(Y'_{p(j)})$ is the identity and we can define $b_{p(j)} \colon F(Y'_{p(j)}) \longrightarrow E$ to be $a_j$.
    \item Otherwise, $b_{p(j)} \colon F(Y'_{p(j)}) \longrightarrow E$ is induced by the universal property of the pushout, as depicted in the following diagram in $\mathsf{E}$
    \[\begin{tikzcd}
	{F(X_i)} & {F(Y_{f(i)})} \\
	{F(Y_{g(i)})} & {F(Y'_{p(j)})} \\
	&& E
	\arrow["{F(\phi_i)}", from=1-1, to=1-2]
	\arrow["{F(\gamma_i)}"', from=1-1, to=2-1]
	\arrow["{F(\pi_{f(i)})}", from=1-2, to=2-2]
	\arrow["{a_{f(i)}}", bend left, from=1-2, to=3-3]
	\arrow["{F(\pi_{g(i)})}"', from=2-1, to=2-2]
	\arrow["{a_{g(i)}}"', bend right, from=2-1, to=3-3]
	\arrow["{b_{p(j)}}"{description}, dashed, from=2-2, to=3-3]
    \end{tikzcd}\]
    (by hypotheses, $F$ preserves pushouts).
\end{enumerate}
The universal property of the coproduct then gives $b \colon \hat{F}(Y'_{J'}) \longrightarrow E$ with the desired property. After the above simplification justifying the limitation to the case $f=g$, with a similar reasoning (this time using that $F$ preserves coequalizers) we get the claim.
\end{proof}

A particular case of \cref{prop:pres} is the following.
\begin{corollary}\label{cor:presh}
If $\mathsf{C}$ has finite connected colimits, then $\Lex(\Fam_\fin(\mathsf{C})^\op,\Set)$ is equivalent to the full subcategory of $\PSh(\mathsf{C})$ consisting of those presheaves $\mathsf{C}^\op \longrightarrow \Set$ that preserve finite connected limits.
\end{corollary}
\begin{proof}
The claim follows by instantiating \cref{prop:pres} with $\mathsf{E} = \Set^\op$.
\end{proof}

At the beginning of this section, we observed that $\Fam_\fin(\Setfs)$ is equivalent to $\SFr_\fin$ (equivalence relations are in one-to-one correspondence with partitions, hence any equivalence relation on a finite set is uniquely determined by a finite family of finite sets; functions between them sending equivalence classes onto equivalence classes are then families of surjective functions among the finite sets of the partitions). Thus we obtain the following.
%silvio 28/10/25 tolti riferimenti ai frame localmente finiti per evitare riferimenti superflui al paragrafo 3
%Moreover, $\SFr_\lf$ is equivalent to $\Lex(\SFr_\fin^\op,\Set)$.
\begin{corollary}\label{cor:catchar}
%$\SFr_\lf$ 
The category $\Lex(\SFr_\fin^\op,\Set)$ is equivalent to the full subcategory of $\PSh(\Setfs)$ consisting of those presheaves $(\Setfs)^\op \longrightarrow \Set$ that preserve finite connected limits.
\end{corollary}
\begin{proof}
We have that $\Lex(\SFr_\fin^\op,\Set)\simeq \Lex(\Fam_\fin(\Setfs)^\op,\Set)$, so that the result follows immediately by 
\cref{cor:presh} if we prove that $\Setfs$ has coequalizers and pushouts. 
%We prove that $\Setfs$ has coequalizers and pushouts, so that the claim follows from the application of \cref{cor:presh} to $\mathsf{C} = \Setfs$. 
To see this, let us identify $\Setfs$ with its image in $\SFr_\fin$, via the full embedding $\Setfs \longrightarrow \SFr_\fin$ sending $\underline{n}$ to the cluster $(\underline{n},\underline{n} \times \underline{n})$. Below we refer to the   computation of  finite colimits in $\SFr_\fin$ 
as described in~\cite[Propositions 3.8 and 4.2]{our_APAL} 
and in the proof of \cref{thm:fromAPAL} above: from such a description it follows that $\Setfs$ is closed under computation of pushouts and coequalizers taken inside $\SFr_\fin$. In detail:
the $\SFr_\fin$-coequalizer 
\[\begin{tikzcd}
	{(\underline{m},\underline{m} \times \underline{m})} & {(\underline{n},\underline{n} \times \underline{n})} & {(W,R)}
	\arrow[shift left, from=1-1, to=1-2]
	\arrow[shift right, from=1-1, to=1-2]
	\arrow[from=1-2, to=1-3]
\end{tikzcd}\]
of a given pair of parallel morphisms in $\Setfs$ is a surjective p-morphism; 
%(see \cite[Propositions 3.8 and 4.2]{our_APAL}); 
this forces the equivalence relation $R$ to be total (p-morphisms map equivalence classes onto equivalence classes) and therefore $(W,R) \in \Setfs$. 
The case of pushouts 
%(which are  p-morphic images of disjoint unions in $\SFr_\fin$) 
is handled in a similar way.
%
%This proves that $\Setfs$ has coequalizers. 
%Using the  description of pushouts in $\SFr_\fin$   as p-morphic images of disjoint unions~\cite[Propositions 3.8 and 4.2]{our_APAL},  we also obtain existence of pushouts in $\Setfs$.
\end{proof}

 \section{Symmetric groups actions}\label{sec:actions}

The plan is to recast the categorical description found in
\cref{cor:catchar}
in more combinatorial terms via finite group actions. 
In this section, we fix some terminology and introduce the required ingredients.

An action (\emph{right action}) of a group $G$ over a set $X$ is an algebraic structure having $X$ as support set and as operations a unary operation for every element $f \in G$;
the action is represented as a function
$$a \colon G \times X \lra X$$
subject to the following  axioms:
$$
a(1,x)=x, \qquad a(f,a(g,x))=a(g\circ f,x)
$$
(where $\circ$ is the binary operation of $G$ and $1$ is the unit of $G$). An action is \emph{transitive} iff for every $x,y \in X$ there is $f \in G$ such that $a(f,x)=y$.
The \emph{orbit} of $x\in X$ is the subset $\{a(f,x)\ \vert\ f \in G\}$; every action is transitive if restricted to an orbit and every action is the disjoint union of the restrictions to its orbits.

% (we indicate directly with $f$ the operation corresponding to $f$). In addition, the following conditions must be satisfied (we usually write $fx$ for $f(x)$):
%\begin{align*}\label{eq:1}
%1x = x, \qquad f(gx) = (g \circ f)x.\tag{1}
%\end{align*}
Given an action as above and $x \in X$, we denote by $\Fix(x)$ the subgroup $\graffe{f \in G\ \vert\ a(f,x) = x}$. The action is said to be \emph{faithful} iff $\Fix(x) = \graffe{1}$ for all $x \in X$. Among the faithful transitive actions, we have the \emph{canonical} action of $G$ on itself $c \colon G\times G\lra G$ defined by $c(f,g)=g\circ f$.

It is well known that an action of a group $G$ can be equivalently seen as a presheaf $G^\op \longrightarrow \Set$, where $G$ denotes both the group itself and the category having as objects the singleton $\star$ and as morphisms $\star \longrightarrow \star$ the elements of $G$ (the binary operation of $G$ is composition); natural transformations between presheaves then give the 
appropriate
%natural 
notion of morphisms between $G$-actions.
%Below, we shall indicate with $G$ both the group $G$ and the one-object category $\underline{G}$.

As we said in the previous section, we denote by $\underline{n}$ the finite set $\graffe{1,\dots,n}$; we are especially interested in the actions of the group $S_n$: this is the \emph{$n$-th symmetric group}, i.e.\ the group of permutations over $\underline{n}$. 

We can lift an action of $S_m$ to a presheaf $(\Setfs)^\op \lra \Set$  as follows.

\begin{definition}\label{def:can}
A presheaf $L \colon (\Setfs)^\op \longrightarrow \Set$ is said to be a \emph{canonical lifting of} a given action $a \colon S_m \times X \lra X$, if it satisfies the following conditions:
\begin{description}
\item[{\rm (i)}] there is a bijection $\eta \colon X \longrightarrow L_m$; %(in what follows, we identify the elements of $L_m$ with the elements of $X$, and simply write $L_m = X$);
\item[{\rm (ii)}] for every $n$,
%matteo $n \geq m$,
the functions $L(q)$'s, varying $q$ among the surjections
$\underline{n} \lra \underline{m}$, cover $L_n$ (in the sense that the union of their images is $L_n$);\footnote{Observe that, if $n < m$, then $L_n$ is the empty set, the set of surjections $\underline{n} \longrightarrow \underline{m}$ being empty.}
\item[{\rm (iii)}] for every 
$q_1,q_2 \colon \un \pair \um$ and $x_1,x_2 \in X$, $L(q_1)(\eta(x_1))=L(q_2)(\eta(x_2))$ holds if and only if there is $\sigma \in S_m$ such that $a(\sigma,x_1)=x_2$ and $\sigma \circ q_2=q_1$.
\end{description}
\end{definition}

To show that the above properties characterize, up to isomorphism, 
%exactly one (if a canonical lifting exists) presheaf $L$.
canonical liftings we show that canonical liftings are nothing but Kan extensions. Recall~\cite[Ch.10]{CWM} that given two functors 
\[\begin{tikzcd}
	{\mathbf{S}} & {\mathsf{C}_1} & {\mathsf{C}_2}
	\arrow["a"', from=1-2, to=1-1]
	\arrow["\iota", from=1-2, to=1-3]
\end{tikzcd}\]
a \emph{left Kan extension} of $a$ along $\iota$ is a pair $(L,\eta)$ given by a functor 
$L \colon \mathsf{C_2}\longrightarrow \mathbf{S}$ and a natural transformation $\eta \colon a \Longrightarrow L\circ \iota$ satisfying the following universal property:
\begin{center}
\say{For every other pair $Y \colon \mathsf{C_2} \longrightarrow \mathbf{S}$ and $\mu \colon a \Longrightarrow Y \circ \iota$, there is a unique natural transformation $\xi \colon L \Longrightarrow Y$ such that $(\xi\cdot\iota)\circ \eta=\mu$.}
\end{center}
\begin{proposition}\label{prop:can1}
If $L \colon (\Setfs)^\op \longrightarrow \Set$ is a canonical lifting of an action $a \colon S_m \times X \lra X$, then $L$ is the left Kan extension of the action $a$ (seen as a presheaf defined on the one-object category $S_m$) along the obvious full and faithful embedding $\iota_m \colon S_m \longrightarrow \Setfs$.
\end{proposition}
\begin{proof}
Consider the function $\eta \colon X \longrightarrow L_m$ in (i) of \cref{def:can}. For each $\sigma \in S_m$ and $x \in X$, we have (by the right-to-left implication of (iii))
\begin{align*}
L(\sigma)(\eta(x)) = L(1)(\eta(a(\sigma,x))) = \eta(a(\sigma,x))
\end{align*}
hence $\eta$ defines a natural isomorphism $\eta \colon a \Longrightarrow L \circ \iota_m$.

To prove the universal property for the pair $(L,\eta)$, consider another presheaf $Y \colon (\Setfs)^\op \longrightarrow \Set$, and a natural transformation $\mu \colon a \Longrightarrow Y \circ \iota_m$ (i.e.\ a function $\mu \colon X \longrightarrow Y_m$ such that $Y(\sigma)(\mu(x)) = \mu(a(\sigma,x))$ for each $\sigma \in S_m$ and $x \in X$). We want to define a natural transformation $\xi \colon L \Longrightarrow Y$ such that $\xi_m \circ \eta = \mu$, and prove that such $\xi$ is unique. To do so, define the function $\xi_n \colon L_n \longrightarrow Y_n$ as follows: given $z \in L_n$, we can write $z = L(q)(\eta(x))$ for some $q \colon \underline{n} \longrightarrow \underline{m}$ and $x \in X$, by (ii) and (i); set $\xi_n(z) \coloneqq Y(q)(\mu(x))$.
%matteo
For each $n$, the function $\xi_n$ is well defined, by the left-to-right implication of (iii), by naturality of $\mu$ and by functoriality of $Y$.\footnote{Such definition makes sense only for $n \geq m$: if $n < m$, as we observed before, $L_n$ is empty, hence $\xi_n \colon L_n \longrightarrow Y_n$ is the empty function.} To see that the collection of the $\xi_n$'s defines a natural transformation $\xi \colon L \Longrightarrow Y$, consider a morphism $p \colon \underline{k} \longrightarrow \underline{n}$ in $\Setfs$ and $z \in L_n$, say $z = L(q)(\eta(x))$ as above. We have that
\begin{align*}
Y(p)(\xi_n(z)) &= Y(p)(Y(q)(\mu(x))) = Y(q \circ p)(\mu(x)) = \xi_k(L(p)(z)),
\end{align*}
since $L(p)(z) = L(p)(L(q)(\eta(x))) = L(q \circ p)(\eta(x))$. Moreover, for each $x \in X$, we have that $\xi_m(\eta(x)) = \xi_m(L(1)(\eta(x))) = Y(1)(\mu(x)) = \mu(x)$. To conclude, if $\xi' \colon L \Longrightarrow Y$ is another natural transformation such that $\xi'_m \circ \eta = \mu$, then, for $z = L(q)(\eta(x)) \in L_n$, we have
\begin{align*}
\xi'_n(z) = \xi'_n(L(q)(\eta(x))) = Y(q)(\xi'_m(\eta(x))) = Y(q)(\mu(x)) = \xi_n(z)
\end{align*}
showing uniqueness of $\xi$.
\end{proof}

Since $S_m$ is small and $\Set$ is co-complete, the left Kan extension mentioned in 
\cref{prop:can1} exists and is computed as a pointwise colimit~\cite[Ch.10,\S 3]{CWM}; we describe it explicitly in the next Proposition.
%For each action $a \colon S_m \times X \longrightarrow X$, a canonical lifting always exists, as shown by the following Proposition. 

\begin{proposition}\label{prop:can2}
Consider an action $a \colon S_m \times X \lra X$ and define a presheaf $L^a \colon (\Setfs)^\op \longrightarrow \Set$ as follows:
\begin{itemize}
     \item[-] we take as $L^a_n$ the set of pairs $(x,q)$ (for $x \in X$ and surjective $q \colon \un \lra \um$)
     divided by the equivalence relation given by $(x_1,q_1)\approx (x_2,q_2)$ iff 
     there is $\sigma \in S_m$ such that $a(\sigma,x_1)=x_2$ and $\sigma \circ q_2=q_1$;
      \item[-] we set $L^a(p)([x,q]) \coloneqq [x,q\circ p]$, for 
      $\underline{k} \buildrel{p}\over \lra \un \buildrel{q}\over\lra \um$ ($[x,q]$ denotes the $\approx$-equivalence class of $(x,q)$).
\end{itemize}
Then $L^a$ is a canonical lifting of the action $a$.
\end{proposition}
\begin{proof}
Showing that $L^a$ is well-defined and that it defines a functor $(\Setfs)^\op \longrightarrow \Set$ is straightforward.

We prove that the above data satisfy conditions (i), (ii) and (iii) of \cref{def:can}. The function $\eta \colon X \longrightarrow L^a_m$ in (i) can be defined as follows: send $x \in X$ to the equivalence class $[x,1] \in L^a_m$ ($1 \colon \underline{m} \longrightarrow \underline{m}$ is the identity). Injectivity of $\eta$ follows from the fact that $[x_1,1] = [x_2,1]$ means, by definition, that there exists $\sigma \in S_m$ such that $a(\sigma,x_1) = x_2$ and $\sigma \circ 1 = 1$, the latter implying $\sigma = 1$, hence $x_2 = a(1,x_1) = x_1$. Moreover, $\eta$ is surjective, since the elements of $L^a_m$ are of the form $[x,\sigma]$, with $\sigma \in S_m$ (surjections $\underline{m} \longrightarrow \underline{m}$ are bijections), and $[x,\sigma] = [a(\sigma,x),1] = \eta(a(\sigma,x))$. To prove (ii), it is sufficient to observe that, given $[x,q] \in L^a_n$, with $x \in X$ and $q \colon \underline{n} \longrightarrow \underline{m}$, we can write $[x,q] = L^a(q)([x,1]) = L^a(q)(\eta(x))$. The latter equation proves also property (iii) for $L^a$.
\end{proof}

In light of \cref{prop:can1}, since Kan extensions are unique up to isomorphisms, we will refer to the  construction of \cref{prop:can2} above  as \emph{the} canonical lifting of an action $a \colon S_m \times X \longrightarrow X$.

\section{The classified theory}\label{sec:classified}

\cref{cor:catchar} showed that the category of models of $T_\textit{S5}^d$ (the theory dual to the theory of \textit{S5}-algebras)
%In Section~\ref{sec:families} we showed that the theory classified by the category of finite $S5$-algebras (i.e. the  theory dual to $T_\textit{S5}$) 
can be equivalently described as as the category of equalizers and pullbacks-preserving  functors $(\Setfs)^\op \longrightarrow \Set$. This category can be presented as the category of models of an essentially algebraic theory $T_0$, by rewriting functoriality and the relevant preservation conditions as axioms in the essentially algebraic fragment of first order logic, along the lines explained in \cref{sec:ea}. However, little information is gained in this way: $T_0$ is just a bit more simple than the logically equivalent theory $T_\textit{S5}^d$
(the latter being obtained by logical transcription of  lex functoriality conditions for contravariant functors from $\SFr_\fin^\op$ to $\Set$).
In order to obtain a more transparent axiomatization, we shall shift to first-order classical logic, in the sense that we shall introduce a \emph{classical}\footnote{ We shall need classical negation to express conditions like $\Fix(x)=\{1\}$.}  first order theory $T_2$ whose axioms 
are \emph{logically equivalent to $T_0$}, so that 
 the models of $T_2$ 
 %(and homomorphisms) 
 will coincide with the models of $T_0$ and of $T_\textit{S5}^d$.
%\footnote{By `homomorphism' we mean the usual morphism in the standard algebraic sense (= well-typed families of functions preserving the operations), in contraposition to the `elementary morphism' (aka elementary embedding) which is the natural notion of a morphism in full first-order logic. Notice that  morphisms preserving truth of the formulae in the essentially algebraic fragment of the logic are the same as algebraic morphisms, so there is no need to use a special name for them to avoid confusion.}

Our strategy comes from the observation that a model $M$ of $T_0$ is a multi-sorted structure consisting of unary maps $M(q) \colon M_n \lra M_m$, varying $q$ among the surjective functions $q \colon \um \lra \un$. When $n=m$, the functoriality conditions show that the $M(q)$ are the actions of the symmetric groups $S_n$; we are going to show, in particular, that such actions \emph{determine uniquely} the structure of the whole model $M$, by using the material introduced in \cref{sec:actions}.

We first fix our working first-order language $\mathcal L$: we have a multi-sorted language (the sorts are the positive natural numbers) with a unary function of type $m\longleftarrow n$ for  every surjective function $f:\um\lra \un$. We use the same letter $f$ to denote both the function $f$ and the corresponding function symbol in the language $\mathcal{L}$ (but recall that domain/codomain get reversed
when we pass from $f:\um\lra \un$ to $f:m\longleftarrow n$).
%to the source/target sorts of  function symbols). 
We reserve the letters $\sigma, \tau, \dots$  to permutations, whereas the letters $f,g,\dots$ are used for generic surjective functions (as such they denote functions that may or may not be bijective). 

The initial theory is the theory $T_0$, which is obtained simply by translating into $\mathcal L$ the functoriality conditions and the conditions of preservation of equalizers and pullbacks of $(\Setfs)^\op$.
We define an intermediate $\mathcal L$-theory $T_1$ (whose axioms  turn out to be logically equivalent to the axioms of $T_0$) as follows.

First of all, in order to  maintain the functoriality conditions,  
we include among  axioms of $T_1$ the equalities 
\begin{align}\label{eq:1}
1(x) = x, \qquad f(g(x))=(g \circ f)(x)\tag{1}
\end{align}
for every $n>0$ (here $x$ is a variable of sort $n$ and $1$ is the identity $\underline{n} \longrightarrow \underline{n}$) and for every $\underline{r}\buildrel{f}\over \lra \um\buildrel{g}\over \lra \un$
in $\Setfs$.

Observe that, up until now, our axioms lie in the essentially algebraic fragment; in order to recapture in a transparent way the preservation of equalizers and pullbacks, we step out from such a fragment. First, we introduce the following abbreviation: for a variable $x$ of sort $n$, we denote by $\Fix(x) = 1$ the formula
\begin{align*}
\bigwedge \{\neg (\sigma (x) = x)\ \vert\ \sigma \in S_n \setminus \{1\}\}.
\end{align*}
Then, the axioms of our theory $T_1$ include
\begin{align*}\label{eq:2}
\bigvee \{\exists y\; (f (y) = x \wedge \Fix(y) = 1)\ \vert\ f \colon \underline{n} \longrightarrow \underline{m}\}\tag{2}
\end{align*}
varying $ n>0$, and
\begin{align*}\label{eq:3}
\left(f (x) = g (y) \wedge \Fix(x) = 1 \wedge \Fix(y) = 1\right) 
%\to \\ 
\to \bigvee \{\sigma (y) = x\ \vert\ 
%\sigma \colon \underline{n} \longrightarrow \underline{m} \text{ bijection s.t. } 
\sigma \circ f  = g\}.\tag{3}
\end{align*}
varying $f \colon \underline{k} \longrightarrow \underline{n}$ and $g \colon \underline{k} \longrightarrow \underline{m}$ in $\Setfs$.
Notice that, according to the above notational conventions, the $\sigma$ in
%(it is intended that 
the conclusion of the axiom \cref{eq:3} is supposed to be a bijection; as a consequence, in the case  $m\neq n$, the disjunction 
occurring in \cref{eq:3} is empty and so it represents
the falsum $\bot$.

\begin{remark}\label{rem:collfree}
Consider a $T_1$-model $M$ (let $M_n$ be the domain associated to the sort $n>0$). Axioms \cref{eq:2} say that
if the action of $S_n$ is not faithful at $x \in M_n$ (i.e.\ if $\Fix(x)$ is not the identity subgroup), then 
it is possible to consider $f \colon \underline{n} \longrightarrow \underline{m}$ and $y \in M_m$ such that the action of $S_m$ is faithful at $y$, and $x = M(f)(y)$
(we use $M(f)$ for the interpretation of the function symbol $f$ in $M$). Moreover, by axioms \cref{eq:3}, such $f$ and $y$ are essentially unique, in the following sense: if $f' \colon \underline{n} \longrightarrow \underline{m'}$ and $y' \in M_{m'}$ satisfy the same property, then 
$m=m'$ and we can find a bijection $\sigma \colon \underline{m} \longrightarrow \underline{m}$ such that $\sigma \circ f = f'$ and $M(\sigma)(y') = y$. Notice that we have exactly $m!$ possible choices for our $f$ and $y$, matching the cardinality of the orbit of $y$ under the faithful action of $S_m$.
\end{remark}

\begin{comment}
\begin{remark}\label{rem:eq}
If the action of $S_n$ is free at $x \in X_n$, then the action is faithful at $x$. Consider $x \in X_n$ and assume that $f x = g x$ for some $f, g \colon \underline{k} \longrightarrow \underline{n}$. By (3), we can find some $q \colon \underline{n} \longrightarrow \underline{m}$ such that $q \circ f = q \circ g$ and $y \in X_m$ such that $x = f y$. If $\Fix(x) = \{1\}$, then $q$ must be bijective by \cref{rem:nodecr}, hence $f = g$.
\end{remark}
\end{comment}

\begin{remark}
Our functions symbols $f \colon m \longrightarrow n$, 
%for $f \colon \underline{n} \longrightarrow \underline{m}$ surjective, 
turn out to  be interpreted as injective functions, i.e. the formula
\begin{align*}\label{eq:4}
f (x) = f (y) \to x = y.\tag{4}
\end{align*}
is true in any $T_1$-model $M$ 
as a consequence of axioms \cref{eq:2} and \cref{eq:3}.
This is seen as follows. Assume
$M\models f(x)=f(y)$; by \cref{eq:2} we have  
$x=M(p)(x')$ and $y=M(q)(y')$ with $M\models\Fix(x')=1$ and $M\models \Fix(y')=1$. From 
$M\models p(f(x'))=q(f(y'))$,
we get, by \cref{eq:3}, that there is $\sigma$ such that 
 $M(\sigma)(y')=x'$ and  $\sigma\circ p\circ f=q\circ f$. Since $f$ is surjective, we obtain $\sigma\circ p=q$. Hence $x=M(p)(x')=M(p)(M(\sigma)(y'))=M(\sigma \circ p)(y')=M(q)(y')=y$.
\end{remark}

\begin{proposition}
The category of
%silvio 31/10
%full subcategory of 
%equalizers and pullbacks
finite connected limits preserving functors $(\Setfs)^\op \longrightarrow \Set$ is equivalent to the category of $T_1$-models and homomorphisms.
\end{proposition}
\begin{proof}
The equivalence between the category of the $\mathcal{L}$-structures which satisfy axioms \cref{eq:1} (and homomorphisms between them) and the category of presheaves $\PSh(\Setfs)$ is straightforward. To prove the claim, we need to show that a functor $X \colon (\Setfs)^\op \longrightarrow \Set$ preserves equalizers and pullbacks if and only if it satisfies (seen as a $\mathcal{L}$-structure) axioms \cref{eq:2,eq:3}.

Assume that $X$ preserves equalizers and pullbacks. 
%\cref{eq:2} follows from the fact that any $f \colon \underline{n} \longrightarrow \underline{m}$ in $\Setfs$ is the coequalizer of some diagram, and that, by hypothesis, such a diagram is mapped by $X$ to a diagram of equalizer in $\Set$, implying the injectivity of $X(f) \colon X_m \longrightarrow X_n$. 
To see \cref{eq:2}, consider $x \in X_n$: if $\Fix(x) = \{1\}$, then there is nothing to prove; otherwise, the existence of $\sigma \neq 1 \in \Fix(x)$ allows us to consider $f$ as the $\Setfs$-coequalizer of $\sigma$ and $1$, so that $y$ is obtained by the preservation of equalizers; using an inductive argument, we can assume that $\Fix(y) = \{1\}$. Finally, consider $f, g$, $x \in X_n$ and $y \in X_m$ as in the premises of \cref{eq:3} and take the pushout of $f$ and $g$ in $\Setfs$, say given by $f' \colon \underline{n} \longrightarrow \underline{k'}$ and $g' \colon \underline{m} \longrightarrow \underline{k'}$; preservation of pullbacks guarantees the existence of some $z \in X_{k'}$ such that $X(f')(z) = x$ and $X(g')(z) = y$. Observe that $f'$ must be bijective (the same holds for $g'$). Otherwise, there is a bijection $\tau \colon \underline{n} \longrightarrow \underline{n}$, different from the identity, acting transitively on the fibers of $f'$. For such $\tau$, we then would have $f' \circ \tau = f'$, hence
$$X(\tau)(x) = X(\tau)(X(f')(z)) = X(f' \circ \tau)(z) = X(f')(z) = x.$$
By $\Fix(x) = \{1\}$, we deduce $\tau = 1$, against our assumptions. To conclude, $\sigma \coloneqq (g')^{-1}f'$ gives the bijection required by the consequent of \cref{eq:3}.

Vice versa, consider a presheaf $X \in \PSh(\Setfs)$ satisfying \cref{eq:2,eq:3}
(hence also \cref{eq:4}). We want to prove that $X$ sends coequalizers and pushouts in $\Setfs$ to, respectively, equalizers and pullbacks in $\Set$. For, consider a diagram of coequalizer
\[\begin{tikzcd}
	{\underline{k}} & {\underline{n}} & {\underline{m}}
	\arrow["f", shift left, from=1-1, to=1-2]
	\arrow["g"', shift right, from=1-1, to=1-2]
	\arrow["q", from=1-2, to=1-3]
\end{tikzcd}\]
in $\Setfs$. We want to show that $X_m$ is isomorphic, via $X(q)$, with the subset given by those $x \in X_n$ such that $X(f)(x) = X(g)(x)$. By \cref{eq:4}, $X(q)$ is injective. To prove its surjectivity, fix any $x \in X_n$ such that $X(f)(x) = X(g)(x)$. By \cref{eq:2}, there exist $p \colon \underline{n} \longrightarrow \underline{l}$ and $y \in X_l$, such that $x = X(p)(y)$ and $\Fix(y) = \{1\}$. By \cref{eq:3}, since $X(p \circ f)(y) = X(p \circ g)(y)$ and $\Fix(y) = \{1\}$, there exists a bijection $\sigma \colon \underline{l} \longrightarrow \underline{l}$ such that $\sigma \circ(p \circ f) = p \circ g$ and $X(\sigma)(y) = y$. The fact that the action of $S_l$ is faithful at $y$ forces $\sigma$ to be the identity and, as a consequence, $p \circ f = p \circ g$. By the universal property of the coequalizer, there exists a unique $r \colon \underline{m} \longrightarrow \underline{l}$ in $\Setfs$ such that $r \circ q = p$. We can then consider $z \coloneqq X(r)(y) \in X_m$, for which $X(q)(z) = X(q)(X(r)(y)) = X(r \circ q)(y) = X(p)(y) = x$. Combining conditions \cref{eq:2,eq:3,eq:4} it is also possible to show that $X$ preserves pullbacks.
\end{proof}

As a consequence,
\begin{corollary}
The category of finite \textit{S5}-algebras, as a lex category, classifies the models of $T_1$ and their homomorphisms as $\mathcal{L}$-structures.
\end{corollary}

Armed by the above result, we introduce our final $\mathcal{L}$-theory $T_2$. The axioms of $T_2$ are just a rewriting of the axioms of $T_1$ up to logical equivalence: 
we first express them in English words, then we supply the translations to first-order formulae  using the information contained in 
%Definition~\ref{def:lifting}(i)-(ii)-(iii) and Proposition~\ref{prop:comparison}.
\cref{def:can}.

%silvio 31/10
We still include formulae~\cref{eq:1} among axioms of $T_2$; as a consequence,
a model $M$ of $T_2$ is in particular a collection of $S_n$-actions $\{a_n \colon S_n\times M_n \lra M_n\}_{n\geq 1}$; as it happens with any action, each $a_n$ is the disjoint union of two actions
$$
a'_n \colon S_n\times M'_n\lra M'_n, \qquad 
a''_n \colon S_n\times M''_n\lra M''_n
$$
where $a'_n$ is the faithful sub-action of $a_n$ (that is, $M'_n$ is the set of the $x \in M_n$ such that $\Fix(x) = \{1\}$). Now the axioms of $T_2$ just say that, 
for every $n$, 
\vskip 2mm
%\centerline{\framebox{\emph{$a_n$ is the disjoint union of the liftings of the $a'_m$ for $m\leq n$.}}
\centerline{\framebox{\emph{$M$ is the disjoint union of the canonical liftings of the $a'_n$ for $n\geq 1$.}}
}
\vskip 2mm \noindent
In first order logic, according to \cref{def:can} and \cref{prop:can1}, this is expressed by the two axioms
\begin{equation}\label{eq:5}
   \bigoplus_{m\leq n}\bigvee \{\exists y\; (f (y) = x \wedge \Fix(y) = 1)\ \vert\ f \colon \underline{n} \longrightarrow \underline{m}\} \tag{5}
\end{equation}
for $n>0$ (here $\oplus$ is exclusive-or expressing disjoint unions) and 
\begin{equation}\label{eq:6}
q_1(x_1)=q_2(x_2) \to \bigvee \{\sigma(q_1)=q_2 \mid 
{\sigma\in S_m},~\sigma\circ q_2 =q_1\} \tag{6}
\end{equation}
for every 
$q_1,q_2 \colon \un \pair \um$.
%In particular $a'_n$ is (trivially) the lifting of itself.
%Notice also that every action of $S_1$ is faithful because $S_1$ is the one-element group.
\vskip 2mm \noindent
 We can then state our final result.
\begin{theorem}\label{thm:T2}
The category of finite \textit{S5}-algebras, as a lex category, classifies the models of $T_2$ and their homomorphisms as $\mathcal{L}$-structures.
\end{theorem}
\begin{proof}
We check that a $\mathcal{L}$-structure $M$ is a $T_1$-model if and only if it is a $T_2$-model.
As a preliminary observation, notice that both models of $T_1$ and of $T_2$ are closed under disjoint unions and isomorphisms.

Assume first that $M$ is a $T_2$-model. By the above mentioned closure under isomorphisms and disjoint unions, we can assume that $M = L^a$ (see \cref{prop:can2}) is the canonical lifting of some faithful action $a \colon S_m \times X \longrightarrow X$. $L^a$ satisfies \cref{eq:1}, being a presheaf.
%To verify \cref{eq:2}, consider $p \colon \underline{k} \longrightarrow \underline{n}$ in $\Setfs$ and $[x_1,q_1], [x_2,q_2] \in L^a_n$. If $L^a(p)([x_1,q_1]) = L^a(p)([x_2,q_2])$, then $[x_1,q_1 \circ p] = [x_2,q_2 \circ p]$, hence there exists $\sigma \in S_m$ such that $a(\sigma,x_1) = x_2$ and $\sigma \circ (q_2 \circ p) = q_1 \circ p$. By the surjectivity of $p$, the latter implies $\sigma \circ q_2 = q_1$, hence $[x_1,q_1] = [x_2,q_2]$. 
To prove \cref{eq:2,eq:3}, we make the following observation. Given $[x,q] \in L^a_n$, for some $x \in X$ and $q \colon \underline{n} \longrightarrow \underline{m}$, and $\tau \in S_n$, we have that $\tau \in \Fix([x,q])$ if and only if $[x,q \circ \tau] = L^a(\tau)([x,q]) = [x,q]$, if and only if there exists $\sigma \in S_m$ such that $a(\sigma,x) = x$ and $\sigma \circ q = q \circ \tau$. The first equation forces $\sigma = 1$, the action $a$ being faithful. Hence we have $\tau \in \Fix([x,q])$ if and only if $q = q \circ \tau$ and, in particular, $\Fix([x,q]) = \{1\}$ if and only if $q$ is an isomorphism. Moreover, for each $\sigma \in S_m$, we have that $[x,\sigma] = [a(\sigma,x),1]$. This means that the elements of $L^a$ for which the action is faithful are those of level $m$, i.e.\ those of the form $[x,1]$, varying $x \in X$. For \cref{eq:2}, pick $[x,q] \in L^a_n$. We have that $[x,q] = L^a(q)([x,1])$ and $\Fix([x,1]) = \{1\}$. Finally, for \cref{eq:3}, if we are given $[x,1], [y,1] \in L^a_m$ (by our previous observations, they are generic elements of $L^a$ for which the action is faithful) and a pair of morphisms $f \colon \underline{n} \longrightarrow \underline{m}$ and $g \colon \underline{n} \longrightarrow \underline{m}$ such that $[x,f] = L^a(f)([x,1]) = L^a(g)([y,1]) = [y,g]$, then there exists $\sigma \in S_m$ such that $a(\sigma,y) = x$, hence $L^a(\sigma)([y,1]) = [y,\sigma] = [a(\sigma,y),1] = [x,1]$ and $\sigma \circ f = g$. 

Vice versa, assume that $M$ is a $T_1$-model. By \cref{eq:1}, $M$ defines a presheaf $(\Setfs)^\op \longrightarrow \Set$. Regard $M$ as a collection of actions $\{a_n \colon S_n \times M_n\lra M_n\}_{n\geq 1}$. Pick any $m \in \mathbb{N}$ and consider the faithful part $a'_m \colon S_m \times M'_m \longrightarrow M'_m$ of $a_m$. The $\mathcal{L}$-substructure of $M$ generated by (any element of) $a'_m$ is given, at level $n$, by the elements $M(q)(x) \in M_n$, varying $q$ among the surjections $\underline{n} \longrightarrow \underline{m}$ and $x \in M'_m$. That such substructure satisfies (i) (with the identity function as $\eta$) and (ii) of \cref{def:can} follows by definition, while property (iii) follows from \cref{eq:3}. This proves that the substructure of $M$ generated by some $a'_m$ is the canonical lifting of $a'_m$. Moreover, \cref{eq:3} ensures that two such substructures, say generated respectively by $a'_m$ and $a'_k$, are disjoint whenever $m \neq k$. To conclude, it is sufficient to observe that, by axiom \cref{eq:2}, the canonical liftings of the $a'_n$, for $n \geq 1$, cover $M$.
\end{proof}

In view of \cref{prop:can1}, we can paraphrase \cref{thm:T2} by saying that \emph{finite S5-algebras classify Kan extensions of faithful actions of the the symmetric groups}. 

From the information we collected, it follows that the category of models of the dual theory of $T_{\textit{S5}}$ can be alternatively  described as the category having as objects the collections of faithful actions of the symmetric groups and as arrows the natural transformations between their left Kan extensions. This observation allows to revisit and to re-interpret \cref{teo:charindcompl}.
In fact, a faithful action is a disjoint union of canonical actions, which are actions of $S_n$ onto itself; representing a canonical action of $S_n$ as a set of cardinality $n$, it is then possible to represent a whole model of $T_2$ as a family of finite sets, i.e. as a locally finite \textit{S5}-Kripke frame. 

Homomorphisms among models of $T_2$ can be represented as families of surjections (i.e. as $p$-morphisms)
for the  reason explained below. Every homomorphism between $T_2$-models is, up to iso, a family of homomorphisms between Kan extensions of two canonical actions. Recall now the universal property of  Kan extensions and notice that
a canonical action of $S_n$ is freely generated by the identity. Thus every natural transformation between the canonical lifting (= Kan extension) of  the canonical action of $S_n$ into the canonical lifting of the canonical action $c_m$ of $S_m$ is uniquely determined by a morphism of the canonical action of $S_n$ into the canonical lifting of $c_m$ at level $n$, that is 
%is uniquely determined 
by the choice of an element of the kind $[x,q]$, i.e. of the kind $[1, q\circ x^{-1}]$ (recall \cref{def:can} and that $x$ must be a permutation because $c_m$ is canonical, too). Thus the morphism can be uniquely represented as a surjection $\un \lra \um$.

The above considerations explain why models of $T_2$ can be represented as sets endowed with an equivalence relation whose equivalence classes are finite (i.e. as locally finite $\textit{S5}$-Kripke frames) and why homomorphisms among models of $T_2$ can be represented as $p$-morphisms in the sense of Kripke semantics. This is certainly a simple and manageable description, but on the other hand such description  hides the rich algebraic structure implicit in these models.

%\begin{remark}
%It is possible to characterize the theory $T_n$ corresponding to the countably many extensions of \textit{S5} by restricting to the actions of the symmetric groups $S_m$ for $m \leq n$. In such cases, one can use a mono-sorted signature, by introducing partial operations.
%\end{remark}

%\begin{remark}
%It should be possible to characterize models of $T$ as continuous action of some infinite topological group ‘subsuming’ all $S_n$ (however, such characterization would not be elementary).
%\end{remark}
\section{Conclusions}\label{sec:conclusions}

The classification result we  obtained in \cref{sec:classified} 
can be seen as a curious mathematical \emph{divertissement}, connecting two seemingly totally unrelated areas (modal logic and finite groups), but on the other hand it is somewhat surprising, in the sense that it reveals how unexpectedly complex combinatorics can be hidden behind the (indeed very cryptic) built-in formalism of lex categories. In \cref{sec:classified}, in order to make sense of the axiomatization mechanically obtained 
%inside the essentially algebraic fragment of first order logic  
via the general instructions of \cref{sec:ea}, we had to resort to sophisticated concepts (Kan extensions) --- \emph{per se} expressible only in some higher order type theory --- and to convert them back to first-order logic via specific additional work (\cref{prop:can1}).

Concerning applications to logic, 
we believe that the characterization we obtained for the \textit{S5} case could be extended to other modal and non classical logics, because the conceptual framework inherited  from~\cite{our_APAL} (forming the base of the investigations of this paper) seems to apply to a quite general context.
However, on the other hand, it is not clear how the essentially algebraic features arising from the analysis of the duals of profinite algebras could interact (if they can interact at all) with standard problems of interest to modal logicians, like the complexity of satisfiability problems, or the design of suitable (also infinitary) calculi, or again the metatheoretical properties such as definability, interpolation~\cite{our_mai}, 
etc.: this is a question that we leave to future research.

From the point of view of universal algebra, the investigation of dual essentially algebraic theories might be an interesting subject by itself (a good starting point could be the exploration of Post lattice in that perspective). We have some feeling that connections to constraint satisfiability area might also arise, via  the notion of a so-called \emph{polymorphism}~\cite{constraints}. In fact, a polymorphism of a finite algebraic structure $A$ is a homomorphims $A^n\lra A$ from some power of $A$ to $A$: as such, it is a typical definable $n$-ary function in the dual of the equational theory axiomatizing the variety generated by $A$.  

\bibliographystyle{plain}
\bibliography{biblio}

\end{document}